\documentclass[smallextended]{svjour3}       
\smartqed  
\usepackage{amsmath}
\usepackage{graphicx}
\usepackage{mathptmx}      
\usepackage[utf8]{inputenc}
\let\phi=\varphi
\usepackage{xcolor}
\usepackage{url}

\newcommand{\system}[1]{\textsf{#1}}

\newcommand{\provernine}[0]{\system{Prover9}}
\newcommand{\macefour}[0]{\system{Mace4}}
\newcommand{\systemontptp}[0]{\system{SystemOnTPTP}}
\newcommand{\paradox}[0]{\system{Paradox}}

\newcommand{\eprover}[0]{\system{E}}
\newcommand{\vampire}[0]{\system{Vampire}}
\newcommand{\Ax}[1]{\mathrm{#1}}
\newcommand{\implicationsymbol}[0]{\rightarrow}
\newcommand{\implication}[2]{{#1} \implicationsymbol {#2}}
\newcommand{\betweennesssymbol}[0]{Z}
\newcommand{\between}[3]{\betweennesssymbol({#1},{#2},{#3})}
\newcommand{\negationsymbol}[0]{\neg}
\newcommand{\negation}[1]{\negationsymbol{#1}}
\newcommand{\Hnew}[0]{H_{\mathrm{new}}}
\newcommand{\Hold}[0]{H_{\mathrm{old}}}
\newcommand{\Huntington}[0]{\mathrm{H}}
\newcommand{\McPheeAx}[1]{\mathrm{M}_{#1}}
\newcommand{\McPhee}[1]{\Sigma_{#1}}
\newcommand{\Cn}[1]{\mathrm{Cn}({#1})}
\newcommand{\Hnewprime}{\Hnew^{\prime}}
\newcommand{\coloneqq}[2]{\mathrel{\mathop:}=}
\newcommand{\strictbetweennesssymbol}[0]{\mathrm{B}}
\newcommand{\strictbetween}[3]{\strictbetweennesssymbol({#1},{#2},{#3})}
\usepackage{multirow}
\newcommand{\Hprime}[0]{H^{\prime}}
\newcommand{\defbetween}[0]{\betweennesssymbol_{\mathrm{def}}}
\newcommand{\defstrictbetween}[0]{\strictbetweennesssymbol_{\mathrm{def}}}

\begin{document}

\title{Sharpening independence results for Huntington's affine geometry}

\titlerunning{Sharpening Huntington's geometry}
\author{Jesse Alama}

\institute{Jesse Alama\at{}\email{}}
\date{Received: date / Accepted: date}
\maketitle
\begin{abstract}{We improve upon Huntington's affine geometry by showing that his independence proofs can be, in some cases, simplified.  We carry out a systematic investigation of the strict notion of betweenness that Huntington employs (the three arguments are supposed to be distinct) by comparing it to McPhee's three axiom systems for the same intended class of structures, which employs weak betweenness (the arguments are permitted to be equal).  Upon closely inspecting the proof that McPhee's axiom systems are equivalent to Huntington's (subject of course to the definition of weak betweenness in terms of strict and vice versa), one finds surprisingly that McPhee's axiom systems have quite different relations to strict betweenness.
}
\subclass{{\and}
}
\end{abstract}

\section{Introduction}\label{sec:intro}

We improve upon some of the results from Huntington's ordered geometry~\cite{huntington1924new} and show a curious interplay, in this setting, between ``weak'' betweenness (the relation among three points that allows some of its arguments to be equal) and ``strict'' betweenness (in which the arguments of betweenness are supposed to be distinct objects).  Huntington improved upon an earlier ordered geometries by reducing the number of postulates from $12$ to $5$ (and showed that the more parsimonious set $H$ of axioms has the same models as the original set), and that $H$ is completely independent in the sense of Moore~\cite{moore1910introduction}.  With automated theorem provers we can confirm Huntington's result and are able to slightly improve upon it by offering a more perspicuous enumeration of the $32$ structures that witness $H$'s complete independence.  (To establish the complete independence of an axiom system $S$, one must produce for each subset $X$ of $S$ a model of $X$ in which every member of $S \setminus X$ is false).  More precisely, all of Huntington's $32$ examples have cardinality $4$, but for only $21$ subsets of $S$ are four elements actually required.

We also aim to address a curiosity pointed out by Huntington in his investigation of one of the axioms.  The curiosity is that while the axiom is necessary, somehow it plays hardly any role.  What is going on?  To address this question, we view Huntington's geometry through the lens of a definitionally equivalent theory, namely McPhee's~\cite{mcphee1962axioms}.


Section~\ref{sec:axioms} introduces Huntington's axiom system $H$ for linear affine geometry (that is, the affine geometry of a single line) that is in focus.  Section~\ref{sec:independence-models} recalls the definition of complete independence of an axiom system and refines Huntington's result that $H$ is completely independent.  Our main contribution can be found in Section~\ref{sec:detached-axiom}, which expands on a curious phenomenon noticed by Huntington about $H$, namely, that one of the axioms is necessary but, at the same time, unused (in light of this curious fact, Huntington calls the axiom ``detached'').  Section~\ref{sec:conclusion} concludes.


\section{Axioms}\label{sec:axioms}

The following are the five axioms ($\Ax{A}$, $\Ax{B}$, $\Ax{C}$, $\Ax{D}$, and $\Ax{9}$) that are in focus in the present paper.  Universal quantifiers implicitly govern all variables.  There is only one sort; the intended interpretation consists of points on a line (though there is no notion of ``line'').  Adopting the notational convention in Pambuccian's survey~\cite{pambuccian2011axiomatics}, the symbol ``$\strictbetweennesssymbol$'' is the ternary relation of betweenness.  For brevity, we use ``$\delta(X)$'' for a sequence $X$ of variables to stand for the obvious conjunction asserting that the variables in $X$ are all distinct.

\begin{definition}[$H$] The set $H$ of axioms consists of the following five sentences:

\begin{description}
\item[$\Ax{A}$] $\implication{\strictbetween{a}{b}{c}}{\strictbetween{c}{b}{a}}$
\item[$\Ax{B}$] $\implication{\delta(a,b,c)}{(\strictbetween{b}{a}{c} \vee \strictbetween{c}{a}{b} \vee \strictbetween{a}{b}{c} \vee \strictbetween{c}{b}{a} \vee \strictbetween{a}{c}{b} \vee \strictbetween{b}{c}{a})}$
\item[$\Ax{C}$] $\implication{\delta(a,x,y)}{\negation{(\strictbetween{a}{x}{y} \wedge \strictbetween{a}{y}{x})}}$
\item[$\Ax{D}$] $\implication{\strictbetween{a}{b}{c}}{\delta(a,b,c)}$
\item[$\Ax{9}$] $\implication{(\strictbetween{a}{b}{c} \wedge \delta(a,b,c,x))}{(\strictbetween{a}{b}{x} \vee \strictbetween{x}{b}{c})}$
\end{description}

Axiom $\Ax{D}$ tells us that Huntington's notion of betweenness is the strict one, according to which the three arguments of betweenness are intended to be distinct.  By contrast, a more relaxed view of betweenness (e.g., that taken by Tarski~\cite{tarski1959elementary,tarski1999tarskis}) permits the arguments to be equal.  (In Section~\ref{sec:detached-axiom} we shall investigate another theory of the affine geometry of a line in which this approach is taken.)

\end{definition}

(It is not a mistake that axioms $1$--$8$ are missing from the list.  Part of Huntington's paper is to replace all $8$ axioms, introduced in an earlier paper, by a single axiom, naturally called $9$.  He proves that new set of $5$ axioms derives all of $1$--$8$. In this paper the more verbose set of axioms will not be considered.)



\section{Very small independence models}\label{sec:independence-models}

The notion of completely independent set of axioms, due to E.~H.~Moore~\cite{}, is a considerable generalization of the familiar notion of independent axiom.  We now recall the definition.

\begin{definition}
A set $S$ of sentences is said to be \emph{completely independent} if, for all subsets $X$ of $S$, the set $(S \setminus X) \cup \{ \neg \phi \mid \phi \in X \}$ is satisfiable.
\end{definition}

Intuitively, a set $S$ of axioms is completely independent when there is no (Boolean) relationship among the axioms.  One can freely negate some (or even all) of the axioms, while keeping the rest, and still find models.  Evidently, a completely independent set is consistent (take $X = \emptyset$) and independent (negating any particular axiom still results in a satisfiable set of sentences).

When concerned with showing independence of a set $S$ of axioms, by ``independence model'' one normally understands a structure that satisfies all the axioms save one, which the model falsifies.  When dealing with the notion of independence model, we take a wider notion of ``independence model'' because we are interested in falsifying every subset of $S$, not just its singletons.

With this understanding of ``independence model'', it is clear that to demonstrate the complete independence of a set $S$ of axioms, one has to produce $2^{|S|}$ independence models.  Even for small sets $S$, the task is evidently rather tedious.  It is remarkable that in $1924$ Huntington did not shy away from the task of showing the complete independence of his axiom system.  (Perhaps the desire to establish complete independence for a previously published theory of betweenness, which had $12$ axioms, drove him to simplify the theory.  He managed to go from $12$ to $5$ axioms.  With patience, one can produce $2^{5} = 32$ independence models; it is another matter to try to produce $2^{12} = 4096$, certainly in the pre-computer age.  Nonetheless, van der Walle actually did this~\cite{vanderwalle1924complete}.)

When one inspects the $32$ independence models Huntington produced to show the complete independence of $H$, one finds that they all have cardinality $4$.  One can improve upon this result by showing that only in some cases do we really need $4$ points; for the rest, we can get by with fewer (sometimes even $1$ point is enough).  Table~\ref{tab:size-of-models} lists those cases where an independence model exists having cardinality strictly less than $4$.  For the sake of uniformity, every independence model has the domain $\{1\}$, $\{1\,2\}$, or $\{1\,2\,3\}$.  Our presentation of the models mimics Huntington's: we list the true (and only the true) betweenness statements in the compact notation $abc$, which is to be understood as: in the model, $\between{a}{b}{c}$ holds.

%
\begin{table}
\begin{tabular}{ccccc|c|r}
A & B & C & D & 9 & Minimal cardinality & Interpretation of betweenness in $\{1,2,3\}$\\
\hline
$+$ & $+$ & $+$ & $+$ & $+$ & 1 & (no true betweennesses)\\
$+$ & $+$ & $+$ & $-$ & $+$ & 1 & 111\\
$+$ & $+$ & $-$ & $+$ & $+$ & 3 & 123, 132, 231, 321\\
$+$ & $+$ & $-$ & $-$ & $+$ & 3 & (all possible betweennesses)\\
$+$ & $-$ & $+$ & $+$ & $+$ & 3 & (no true betweennesses)\\
$+$ & $-$ & $+$ & $-$ & $+$ & 3 & 121\\
$-$ & $+$ & $+$ & $+$ & $+$ & 3 & 123\\
$-$ & $+$ & $+$ & $-$ & $+$ & 2 & 111, 122\\
$-$ & $+$ & $-$ & $+$ & $+$ & 3 & 123, 213, 231\\
$-$ & $+$ & $-$ & $-$ & $+$ & 3 & 111, 123, 132, 211 \\
$-$ & $-$ & $+$ & $-$ & $+$ & 3 & 111, 211\\
\end{tabular}
\caption{\label{tab:size-of-models}Only those cases where a smaller cardinality independence model exists are presented.}
\end{table}


\section{A detached axiom}\label{sec:detached-axiom}

Examining the role that axiom $\Ax{D}$ plays in showing that the new axiom set $H$ derives the old axiom set $H \cup \{ \Ax{1}, \dots, \Ax{8} \}$, Huntington notices that $\Ax{D}$ was never used in any proof and remarks:

\begin{quotation}
It will be noticed that postulate D plays a peculiar rôle. Although it is strictly independent and therefore cannot be omitted, yet it is not used in proving any of the theorems on deducibility, and it may always be made to hold or fail without affecting the holding or failing of any other postulate. It may therefore be called not only independent but altogether ``detached''.~\cite[p.~$275$]{huntington1924new}
\end{quotation}

{\noindent}``Theorem on deducibility'' here means the class of theorems showing that the new axiom set $H$ is deductively just as strong as the old axiom set (the precise list of whose axioms need not concern us here).  Huntington does not expand on this interesting comment; after mentioning the issue, he moves on to the heart of his paper, which is to show the complete independence of his axiom set.

The result is indeed somewhat curious.  It appears $\Cn{\Hnew \setminus \{ \Ax{D} \}} = \Cn{\Hold \setminus \{ \Ax{D} \}}$ )that is, the old and the new axiom systems have the same class of models, putting aside $\Ax{D}$, and $\Ax{D}$ is an independent axiom of both.  It suggests that what we are really dealing with is a theory of affine linear order in which betweeness of three points need not imply that the points are distinct.  In other words, $\Hnew \setminus \{ \Ax{D} \}$ is consistent with the falsity of $\Ax{D}$.

\begin{definition}[$\Hprime$]
Define the axiom system $\Hprime$ as $H \setminus \{ \Ax{D} \}$.
\end{definition}

How, then, can one distinguish between $H$ and $\Hprime$?  One cannot distinguish the smallest possible models, because both $\Hnew$ and $\Hnewprime$ have $1$-element models.  (No axiom ensures non-triviality; it is consistent with both theories that there are no true betweenness statements at all.)  The theories can be distinguished when one starts adding hypotheses, such as ``There is at least one true betweenness statement''.  Under this hypothesis, the smallest model of $\Hnew$ has cardinality $3$, whereas there is model of $\Hnewprime$ has cardinality $1$.

In the absence of $\Ax{D}$, a curious possibility is revealed.  The fact that $\Ax{D}$ is not usd to show that $\Cn{\Hold} = \Cn{\Hnew}$ vaguely suggests that perhaps the distinctness of the points is not essential, and perhaps one gets the impression that a generalization of some kind of available.  The intention of the betweenness relation in Huntington's axioms is that the arguments of the point are mutally distinct.  One sees this in Huntington's discussion of his axioms, but conretely, it is obviously by inspecting $\Ax{D}$.

To investigate this curious issue, we propose to view Huntington's axioms through the lens of another theory.  By doing so we are following in Huntington's footsteps.  Evidently, the curiousness of $\Ax{D}$ was seen only through a detailed inspection of various theorems of $H$ (more precisely, the derivations from $H$ of the older axioms $\Ax{1--8}$ that are missing from $H$).  We propose to use McPhee's axiom system for betweenness~\cite{mcphee1962axioms}.  McPhee developed the theory of linear affine geometry using the single undefined notion of ``weak'' betweenness $\betweennesssymbol$, which, unlike Huntington's $\strictbetweennesssymbol$, permits its arguments to be equal.  We now list McPhee's axioms seven axioms and the three axiom systems that he defines in terms of them:

\begin{description}
\item[$\McPheeAx{1}$] $\exists z \forall a, b \left [ \between{a}{b}{z} \vee \between{a}{z}{b} \vee \between{b}{z}{a} \vee \between{b}{a}{z} \vee \between{z}{a}{b} \vee \between{z}{b}{a} \right ]$
\item[$\McPheeAx{2}$] $\implication{(\between{b}{a}{c} \wedge \between{c}{d}{a})}{\between{d}{a}{b}}$
\item[$\McPheeAx{3}$] $\implication{(\between{b}{a}{c} \wedge \between{d}{b}{a})}{(\between{c}{a}{d} \vee a = b)}$
\item[$\McPheeAx{4}$] $\implication{(\between{b}{a}{c} \wedge \between{c}{a}{d} \wedge \between{d}{a}{b})}{(a = b \vee a = c \vee a = d)}$
\item[$\McPheeAx{5}$] $\between{a}{b}{c} \vee \between{b}{c}{a} \vee \between{b}{a}{c}$
\item[$\McPheeAx{6}$] $\between{a}{b}{c} \vee \between{a}{c}{b} \vee \between{b}{c}{a} \vee \between{b}{a}{c} \vee \between{c}{a}{b} \vee \between{c}{b}{a}$
\item[$\McPheeAx{7}$] $\implication{(\between{c}{a}{d} \wedge \between{c}{b}{d} \wedge \between{a}{c}{b} \wedge \between{a}{d}{b})}{(a = c \vee b = c \vee a = d \vee b = d)}$
\end{description}

From these $7$ axiom McPhee puts forward three distinct axiom systems for consideration:

\begin{description}
\item[$\McPhee{1}$] $\{ \McPheeAx{1}, \McPheeAx{2}, \McPheeAx{3}, \McPheeAx{4} \}$
\item[$\McPhee{2}$] $\{ \McPheeAx{3}, \McPheeAx{4}, \McPheeAx{5} \}$, and
\item[$\McPhee{3}$] $\{ \McPheeAx{2}, \McPheeAx{6}, \McPheeAx{7} \}$.
\end{description}

McPhee proves:

\begin{theorem}
$\Cn{\McPhee{1}} = \Cn{\McPhee{2}} = \Cn{\McPhee{3}}$.
\end{theorem}

In other words, McPhee's three axiom systems are equivalent.

One can also show that McPhee's three axiom systems and Huntington's axiom equivalent to each other in the sense that, starting from one, one can define the sole primitive notion of the other and prove all its axioms.  Consider the two interdefinitions of $\betweennesssymbol$ and $\strictbetweennesssymbol$:

\begin{description}
\item[$\defbetween$] $\forall x,y,z \left [ \between{x}{y}{z} \leftrightarrow (\strictbetween{x}{y}{z} \vee x = y \vee x = z \vee y = z \right ]$
\item[$\defstrictbetween$] $\forall x,y,z \left [ \strictbetween{x}{y}{z} \leftrightarrow (\between{x}{y}{z} \wedge x \neq y \wedge x \neq z \wedge y \neq z \right ]$
\end{description}

One can prove (we omit the proofs here):

\begin{proposition}\label{mcphee-implies-huntington}
$\Cn{\McPhee{1} \cup \{ \defstrictbetween \}} = \Cn{\Huntington}$ (and likewise for $\McPhee{2}$ and $\McPhee{3}$).
\end{proposition}

\begin{proposition}\label{huntington-implies-mcphee}
$\Cn{\Huntington \cup \{ \defbetween \}} = \Cn{\McPhee{1}}$ (and likewise for $\McPhee{2}$ and $\McPhee{3}$).
\end{proposition}

To better understand the significance of the ``detached'' axiom $\Ax{D}$, we can look, as Huntington did, into the details of the equivalence of $H$ and McPhee's axioms.  Our question is: What is the analogue of $\Ax{D}$ in McPhee's axiom systems?  Recapitulating Propositions~\ref{mcphee-implies-huntington} and~\ref{huntington-implies-mcphee} with the help of a theorem prover, one proceeds in three steps:

\begin{enumerate}
\item Postulate Huntington's axioms,
\item Postulate $\defbetween$, and
\item Derive (each of) McPhee's axiom systems, inspecting, for each proof, whether $\Ax{D}$ is semantically needed (its absence leads to countersatisfiability).
\end{enumerate}

Table~\ref{tab:huntington-implies-mcphee-d} summarizes the results of this experiment.

\begin{table}
\centering
\begin{tabular}{c|c|c}
Axiom System & Axiom to derive & $\Ax{D}$ needed?\\
\hline
\multirow{4}{*}{$\McPhee{1}$} & $\McPheeAx{1}$ & no\\
 & $\McPheeAx{2}$ & yes\\
 & $\McPheeAx{3}$ & no\\
 & $\McPheeAx{4}$ & yes\\
\hline
\multirow{3}{*}{$\McPhee{2}$} & $\McPheeAx{3}$ & no\\
& $\McPheeAx{4}$ & no\\
& $\McPheeAx{5}$ & yes\\
\hline
\multirow{3}{*}{$\McPhee{3}$} & $\McPheeAx{2}$ & no\\
& $\McPheeAx{6}$ & yes\\
& $\McPheeAx{7}$ & yes\\
\end{tabular}
\caption{\label{tab:huntington-implies-mcphee-d}The role played by axiom $\Ax{D}$ in deriving McPhee's axioms}
\end{table}


Among McPhee's three axiom systems, we find that only $\McPhee{2}$ requires $\Ax{D}$ in the justification of exactly one of its axioms.  For $\McPhee{1}$ and $\McPhee{2}$ we find, interestingly, that the influence of (or need for) $\Ax{D}$ is ``spread out'' or ``smeared'' among multiple axioms.

Going the other way around, let us reconsider the proof of Proposition~\ref{mcphee-implies-huntington} that McPhee's three axiom systems imply Huntington's (with the help of $\defstrictbetween$).  The question, more precisely, is: Starting with one of McPhee's axiom systems, what is necessary to prove $\Ax{D}$?  Table~\ref{tab:mcphee-implies-huntington-d} summarizes the results.

\begin{table}
\centering
\begin{tabular}{c|c|c}
Axiom System & Axiom & Needed to prove $\Ax{D}$\\
\hline
\multirow{4}{*}{$\McPhee{1}$} & $\McPheeAx{1}$ & no\\
& $\McPheeAx{2}$ & no\\
& $\McPheeAx{3}$ & no\\
& $\McPheeAx{4}$ & yes\\
\hline
\multirow{3}{*}{$\McPhee{2}$} & $\McPheeAx{3}$ & no\\
& $\McPheeAx{4}$ & yes\\
& $\McPheeAx{5}$ & no\\
\hline
\multirow{3}{*}{$\McPhee{3}$} & $\McPheeAx{2}$ & yes\\
& $\McPheeAx{6}$ & yes\\
& $\McPheeAx{7}$ & yes\\
\end{tabular}
\caption{\label{tab:mcphee-implies-huntington-d}Necessary axioms for deriving Huntington's axiom $\Ax{D}$ from McPhee's axioms}
\end{table}


The case of $\McPhee{3}$ is quite curious.  There, all axioms are needed (together of course with $\defstrictbetween$) to derive $\Ax{D}$.  That is, if any of $\McPhee{3}$'s axioms is taken away, $\Ax{D}$ cannot be proved.  Thus, it seems intuitively clear that none of $\McPhee{3}$'s axioms can be taken as the analogue of $\Ax{D}$.  As we saw earlier, $\McPhee{2}$ does seem to give a satisfactory answer: only one of its axioms is needed to derive $\Ax{D}$.  But curiously, the ``analogue'' of $\Ax{D}$ here differs from the ``analogue'' we found earlier when we tried to derive $\McPhee{2}$ from $H$ (second bundle of rows in Table~\ref{tab:huntington-implies-mcphee-d})!


\section{Conclusion and future work}\label{sec:conclusion}

We have investigated the extent to which Huntington's axiom $\Ax{D}$, which enforces the distinctness of the three arguments of the betweenness relation, is ``detached'' (as Huntington puts it) by looking at it through the lens of three other theories for the affine structure of a line due to McPhee.  Axiom $\Ax{D}$ may thus be detached, but we hope to have given more insight into how multifaceted its ``grip'' truly is.

A number of further problems suggest themselves.  Pambuccian has drawn attention~\cite{pambuccian2004early} to Huntington's paper (among others) as an early example of ``resource-consciousness'' in mathematics.  In many of Huntington's proofs, it is explicitly indicated not only what premises are used in the proof, but how many times any particular premise is used.  Although this information was ignored for the results reported here, this information is quite interesting and deserves to be further studied as a suite of genuinely mathematical (as opposed to purely logical) cases where close attention is paid to how an axiom is used (and not just whether it is used at all).


Not all geometry is presented in such an explicit, axiomatic style as we find in Huntington's work (and the work of many others).  Indeed, it must be said that only a corner of geometry is given as collections of first-order axiomatic systems.  Owing to expressive limitations of first-order logic, this is probably the best one can hope for.  Nor does a strict axiomatic approach seem to be suitable for all cases, or for all geometers.  Still, when a system of geometry \emph{is} presented as a finitely axiomatized first-order theory, a world of possibilities is opened up owing the the applicability of techniques coming from theorem proving and artificial intelligence.  Research tasks that engaged researchers in the past, such as determining (complete) independence, become ``mechanizable'' (though of course there is never a guarantee that, in general, such problems are readily solvable).  How many problems about geometry were never posed by great mathematical minds owing to the combinatorial horror of actually tackling them?  How many productive mathematical concepts have been put forward, only to flounder (like Moore's notion of completely independent axiom system) owing to their inherent tediousness? The price of exploration now is much lower than it was in Huntington's day.  Let us hope that research in geometry and theorem proving continues to be a fruitful liaison.


\begin{acknowledgements}
The results reported here were achieved thanks to the hard work that goes into developing and maintaining a number of automated reasoning systems, in particular the {\provernine}/{\macefour}\footnote{\url{http://www.cs.unm.edu/~mccune/prover9/}}suite by W.~McCune~\cite{prover9-mace4}, G.~Sutcliffe's {\systemontptp}~\cite{DBLP:conf/lpar/Sutcliffe10}\footnote{\url{http://www.cs.miami.edu/~tptp/cgi-bin/SystemOnTPTP}}, K.~Claessen and N.~Sörensen's {\paradox} suite~\cite{claessen2003new}\footnote{\url{https://github.com/nick8325/equinox}}, S.~Schulz's {\eprover}~\cite{schulz2002e}\footnote{\url{http://www.eprover.org}}, and A.~Voronkov's {\vampire} system~\cite{Vampire}\footnote{\url{http://www.vprover.org}}.

\end{acknowledgements}

\bibliographystyle{spmpsci}      
\bibliography{geometric-reasoning}   

\end{document}